\documentclass[graybox]{svmult}

\usepackage{type1cm}        % activate if the above 3 fonts are
\usepackage{makeidx}         % allows index generation
\usepackage{graphicx}        % standard LaTeX graphics tool
\usepackage{multicol}        % used for the two-column index
\usepackage[bottom]{footmisc}% places footnotes at page bottom

\usepackage%[final]%
{changes}
\definechangesauthor[name=Christian, color=blue!80!black]{ce}
\definechangesauthor[name=Sorin, color=orange!80!black]{sp}
\definechangesauthor[name=Thomas, color=green!80!black]{tw}

\usepackage{newtxtext}       % 
\usepackage{newtxmath}       % selects Times Roman as basic font

\makeindex             % used for the subject index
\usepackage{algorithm}
\usepackage{algpseudocode}
\usepackage{siunitx}

\newcommand{\bsig}{\boldsymbol{\sigma}}
\newcommand{\be}{\mathbf{e}}
\newcommand{\bP}{\mathbf{P}}
\newcommand{\blambda}{\boldsymbol{\Lambda}}

\let\epsilon\varepsilon
\let\phi\varphi
\def\eps{\varepsilon}

\begin{document}

\title*{Dynamic and weighted stabilizations of the $L$-scheme applied to {a} phase-field {model for} fracture {propagation}}
%\title*{Dynamic and weighted stabilizations of the $L$-scheme applied to phase-field fracture}
\titlerunning{Dynamic $L$ scheme for phase-field fracture} 
\author{Christian Engwer, Iuliu Sorin Pop, Thomas Wick}
\institute{Christian Engwer \at Institut f\"ur Numerische und Angewandte Mathematik,
  %Fachbereich Mathematik und Informatik der
  Universit\"at M\"unster,
Einsteinstrasse 62,
48149 M\"unster, Germany, \email{christian.engwer@uni-muenster.de}
\and Iuliu Sorin Pop \at 
Universiteit Hasselt, Faculty of Sciences, Agoralaan Gebouw D - B-3590 Diepenbeek, Belgium, \email{sorin.pop@uhasselt.be}
\and Thomas Wick \at Leibniz Universit\"at Hannover, Institut f\"ur Angewandte Mathematik,
Welfengarten 1, 30167 Hannover, Germany,
\email{thomas.wick@ifam.uni-hannover.de}
}
\maketitle

%%%%%%% warum zwei mal?
% \abstract*{This contribution is a result from immediate discussions 
% after the oral ENUMATH presentation on the $L$-scheme with constant parameters.
% We consider a phase-field fracture propagation problem, which consists 
% of two (nonlinear) coupled PDEs. The first equation describes 
% the displacement evolution and the second equation is a smoothed 
% indicator variable. The latter one is also subject 
% to a variational inequality constraint in time. We propose 
% an iterative scheme, the so-called $L$-scheme, with a dynamic 
% update of the stabilization parameters. Our algorithmic 
% developements are substantiated with two numerical tests.
% Therein, we observe a significant reduction of iteration numbers 
% in comparison to the constant choice of the stabilization parameters.
% }

\abstract{
We consider a phase-field fracture propagation model, which consists 
of two (nonlinear) coupled partial differential equations. 
The first equation describes 
the displacement evolution, and the second is a smoothed 
indicator variable, describing the crack position. We propose 
an iterative scheme, the so-called $L$-scheme, with a dynamic 
update of the stabilization parameters during the iterations. 
Our algorithmic improvements are substantiated with two numerical tests.
The dynamic adjustments of the stabilization parameters lead
to a significant reduction of iteration numbers 
in comparison to constant stabilization values.
}

%%%%%%%%%%%%%%%%%%%%%%%%%%%%%%%%%%%%%%%%%%%%%%%%%%%%%%%%%%
\section{Introduction}
This work is an extension of \cite{BrWiBeNoRa19} in which 
an L-type iterative scheme 
(see \cite{list2016study, MR2079503}) with stabilizing parameters for solving 
phase-field fracture problems was proposed. In \cite{BrWiBeNoRa19}, 
the stabilization parameters were chosen as constants throughout an entire 
computation. With these choices,
the convergence of the scheme has been proven rigorously.
The resulting 
approach performs well in the sense that an unlimited number
of iterations compared to a truncated scheme yields 
the same numerical solution. The results were validated by investigating
the load-displacements curves. Moreover,
the robustness of the scheme w.r.t. spatial mesh refinement was shown. 
Nonetheless, the iteration numbers (for an unlimited number of iterations) 
remained high. 

In this work, we propose and compare two extensions 
of the aforementioned scheme. First,
{we} update the $L$ scheme parameters 
{dynamically}. Second, {we use} 
an adaptive weight depending 
on the fracture location inside the domain. For the latter idea,
we use the phase-field variable to weight $L$ locally.

The outline of this work is as follows: 
In Section \ref{sec_model} the model is stated whereas
Section \ref{sec_L_scheme} presents the  
dynamic choice of the stablization parameters.
In Section \ref{sec_tests}, we present
two numerical tests to study the performance 
of the proposed scheme.

%%%%%%%%%%%%%%%%%%%%%%%%%%%%%%%%%%%%%%%%%%%%%%%%%%%%%%%%%%
\section{The phase-field fracture model}
\label{sec_model}
We consider the crack propagation model proposed in \cite{BrWiBeNoRa19}. 
Let $\Omega \subset \mathbb{R}^d$ be a $d$-dimensional, polygonal and 
bounded domain and $T > 0$ a maximal time. We use common notations, 
and in particular $W^{1, \infty}(\Omega)$ denotes 
the space of functions on $\Omega$ having essentially bounded weak 
derivatives in any direction, while $H^1_0(\Omega)$ contains the functions 
having square integrable weak derivatives, and vanishing at the boundary 
of $\Omega$ (in the sense of traces). For the ease of writing, we let 
$V := H^1_0(\Omega)^d$ and $W := W^{1, \infty}(\Omega)$. 
For (almost every) location $x \in \Omega$ and time $t \in (0, T]$ 
the vector-valued 
displacements are denoted by $u$. The fracture and its propagation 
within $\Omega$ are modeled with the help of a phase field variable 
$\phi$, which approximates the characteristic function of the intact 
region of $\Omega$. Written in weak form, the fracture propagation 
model in resumes to finding 
$(u(t),\phi(t)) \in V\times 
W :=  (H_0^1(\Omega))^d \times W^{1,\infty}(\Omega)$ 
such that for $t \in (0,T]$:
\begin{align}
\nonumber &\textnormal{\textbullet \ \textbf{Step 1}: 
given $(u^{n,i-1}, \phi^{n,i-1})$ find $u^{n,i}$ such that}\\
&a_u(u^{n,i},v) := 
L_u(u^{n,i} - u^{n,i-1},v) + \left( g(\phi^{n,i-1}) \bsig^+(u^{n,i}), \be(v) \right) 
+ \left( \bsig^-(u^{n,i}), \be(v) \right) = 0, \nonumber  \\
&\quad \forall v \in V_h, \label{iter2_mod} \\
\nonumber&\textnormal{\textbullet \ 
\textbf{Step 2}: given $(\phi^{n,i-1}, u^{n,i}, \phi^{n-1})$ find $\phi^{n,i}$ such that} \\
\nonumber &a_{\phi}(\phi^{n,i},\psi) := L_\phi (\phi^{n,i} - \phi^{n,i-1},\psi) + G_c \epsilon(\nabla \phi^{n,i}, \nabla \psi) - \frac{G_c}{\epsilon} (1 - \phi^{n,i}, \psi) \\
&\qquad \qquad+ (1-\kappa)(\phi^{n,i} \bsig^+(u^{n,i}):\be(u^{n,i}),\psi) 
+ ( \Xi + \gamma [\phi^{n,i} - \phi^{n-1}]^+, \psi)= 0, \nonumber \\
&\quad \forall \psi \in W_h. \label{iter1_mod} 
\end{align}
Here, we note that the `time' $t$ appears only the irreversibility constraint
$\partial_t \varphi \leq 0$, yielding an incremental problem and 
which is regularized using an  
augmented Lagrangian penalization $(\Xi + \gamma [\phi^{n,i} - \phi^{n-1}]^+, \psi)$
as proposed in \cite{WheWiWo14}. Here, $\Xi$ is an $L^2(\Omega)$ 
function and $\gamma$ a positive parameter.

Furthermore, in the above, 
$\epsilon$ is a (small) phase-field regularization parameter, 
$G_c > 0$ is the critical elastic energy 
restitution rate, and $0 < \kappa \ll 1$ is 
a regularization parameter used to avoid the degeneracy of the elastic energy. 
The latter is similar to replacing the fracture with a softer material. 
Next, $g(\phi) := (1-\kappa)\phi^2 + \kappa$ is the degradation function, 
and $\be:=\frac 1 2 (\nabla u + \nabla u^T)$ is the 
strain tensor.

The stress tensor in the above is split into a tensile and compressive part, 
\begin{equation*}
\bsig^+ := 2\mu_s \be^+ + \lambda_s [\text{tr}(\be)]^+ I, \quad
\bsig^- := 2\mu_s (\be - \be^+) + \lambda_s \bigl(\text{tr}(\be) -
[\text{tr}(\be)]^+ \bigr) I,
\end{equation*}
where $[\cdot]^+$ stands for the positive cut of the argument. Further,  
$
\be^+ = \bP \blambda^+ \bP^T,
$ 
with $\bP$ being the matrix containing 
the unit eigenvectors corresponding to the 
eigenvalues of the strain tensor $\be$. 
In particular, for $d=2$ one has 
$\bP = [v_1, v_2]$ and 
\[
\blambda^+ := \blambda^+ (u):=
\begin{pmatrix}
[\lambda_1(u)]^+ & 0 \\
0 & [\lambda_2(u)]^+
\end{pmatrix}
. 
\]

%%%%%%%%%%%%%%%%%%%%%%%%%%%%%%%%%%%%%%%%%%%%%%%%%%%%%%%%%%
\section{The $L$-scheme with dynamic updates of the stabilization parameters}
\label{sec_L_scheme}
The iteration \eqref{iter2_mod}-\eqref{iter1_mod} is essentially the scheme proposed in \cite{BrWiBeNoRa19}, 
in which the stabilization parameters $L_u$ and $L_\phi$ 
are taken constant. To improve the convergence behaviour 
of the scheme, we propose a dynamic update of these parameters.

\medskip\noindent\textbf{Dynamic update at each iteration / constant in space:\quad}
The iteration discussed in \cite{BrWiBeNoRa19} uses constant parameters $L_u$ and $L_{\phi}$. 
With this choice, the convergence has been proved rigorously. However, the number of iterations can remain high. High iteration numbers 
for phase-field fracture problems were also reported in \cite{GeLo16,Wi17_SISC}.
To improve the efficiency, we suggest in this work to 
update $L_u$ and $L_{\phi}$ at each 
iteration $i$:
\[
L_{i} = a(i)  L_{i-1}, \qquad \text{where } L_{i} := L_{u,i} = L_{\phi,i}. 
\]
Inspired by 
numerical continuation methods in e.g.  \cite{AllGe90}, one would naturally 
choose a large $L_0$ and $a(i):=a<1$ to obtain a decreasing sequence
$L_0 > L_1 > L_2 > \ldots $, updated until a lower bound $L_-$ is reached.
However, this seems not to be a good choice in phase-field fracture since 
the system does not have a unique solution. 
Consequently, with increasing $i$ the iterations would oscillate in approaching one or another solution, and the algorithm convergence deteriorates. 
For this reason, we propose the other way around: 
the closer the iteration is to 
some solution, the larger the stabilization parameters is chosen, so that the iterations remain close to this solution. 
We choose $a(i):=a>1$, yielding 
$L_0 < L_1 < L_2 < \ldots $ up to a maximal $L_*$. 

\medskip\noindent\textbf{On the specific choice of the parameters:\quad}
A possible choice 
for $a$ is $a(i) := 5^i$ ($i=0,1,2,\ldots$), while $L_0 := 10^{-10}$. This heuristic choice and may be improved by using the solution 
within the iteration procedure, or a-posteriori  error estimates  for the iteration error. 
Moreover, $a(i) := 5^i$
is motivated as follows. Higher values greater than $5$ would 
emphasize too much the stablization. On the other hand, too low values, do not 
lead to any significant enhancement of the convergence behaviour. We substantiate 
these claims by also using $a(i)=10^i$ and $a(i)=20^i$ in our computations.

\medskip\noindent\textbf{{Dynamic update {using the iteration}}:\quad}
An extension of the strategy is to {adapt} the $L$-scheme parameters
in space by using the phase-field variable $\varphi^{n, i-1}$. We still take 
$L_{i} = a L_{i-1}$,
but now $a := a(i, \varphi^{n, i-1})$. Away from the fracture, we have 
$\varphi \approx 1$ and essentially only the elasticity component \eqref{iter1_mod} is being solved. On the other hand, the 
stabilization is important in the fracture region, for which we take 
\[
L_{i} = a(i, \varphi^{n, i-1}) L_{i-1}, \qquad \text{ with } a(i, \varphi^{n, i-1}) := (1 - \varphi^{n, i-1}) a .
\]
Recalling that the fracture is characterised by $\varphi \approx 0$, it becomes clear that the stabilization parameters are acting mainly in the fracture region. 
Finally, to improve further the convergence behaviour of the scheme we adapt $\Xi$ at each iteration. In this case we take
$
\Xi_i = \Xi_{i-1} + \gamma [\phi^{n, i-1}-\phi^{n-1}]^+ . 
$
\begin{algorithm}[thp]
  \label{algo:Lscheme}
  \caption{Dynamic variant of the L-scheme for a phase-field fracture}
\begin{algorithmic}
\State At the loading step $t^n$
\State Choose $\gamma>0$,  $a>1$, as well as $\Xi^0$ and $L_0$. Set $i = 0$.
\Repeat
\State Let $i = i+1$; 
%\State Iterate on $i$ (combined loop)
\State Solve the two problems, namely
\State\;\;\;\;   Solve the (nonlinear) elasticity \eqref{iter2_mod}
\State\;\;\;\;   Solve the nonlinear phase-field \eqref{iter1_mod}
\State Update
%\State\;\;\;\;
~~$L^{i} = a L^{i-1}$
\State Update
%\State\;\;\;\;
~~$\Xi^{i} = \Xi^{i-1} + \gamma [\phi^{n,i} - \phi^{n-1}]^+$
\Until 
% \begin{equation*}
% \label{tol_3}
%  \max(\|u^{n,i} - u^{n,i-1}\|,\|\phi^{n,i} - \phi^{n,i-1}\|,\|\Xi^{i} - \Xi^{i-1}\|) 
%\le \operatorname{TOL}
\State\;\;\;\;   $\max(\|a_u(u^{n,i},v_k)\|,\|a_\phi(\phi^{n,i},\psi_l)\|) \leq \operatorname{TOL},$
% \end{equation*}
\State Set ~~$(u^n,\varphi^n):= (u^{n,i}, \phi^{n,i})$.
\State Increment ~~$t^n \rightarrow t^{n+1}$.
\end{algorithmic}
% \end{Algo}
\end{algorithm}

\medskip\noindent\textbf{The final algorithm:\quad}
The algorithm is based on the iterative procedure for phase-field fracture originally proposed 
in \cite{WheWiWo14}. Therein, the inequality constraint is realized 
by an augmented Lagrangian iteration. Within this loop we update 
the $L$ scheme parameters too. {The resulting is
  sketched in {Algorithm} \ref{algo:Lscheme}, in which $\operatorname{TOL} = 10^{-6}$ is taken, and $L=L_u=L_{\phi}$.
\begin{remark}
For the solution of both nonlinear subproblems \eqref{iter2_mod}
and \eqref{iter1_mod}, we use a monotonicity-based 
Newton method (details see e.g., in \cite{Wi17_SISC}) 
with the tolerance $10^{-8}$.
Inside Newton's method, we solve the linear systems with a direct solver.
\end{remark}

%%%%%%%%%%%%%%%%%%%%%%%%%%%%%%%%%%%%%%%%%%%%%%%%%%%%%%%%%%
\section{Numerical tests}
\label{sec_tests}
We consider two test examples. 
Details for the first test van be found in~\cite{MieWelHof10a}.
The setup of the second test can be found for instance in~\cite{MesBouKhon15}.
Both examples were already computed in~\cite{BrWiBeNoRa19} and 
the results therein 
are compared to the ones obtained here.
The scheme is implemented in a code based on 
the deal.II library~\cite{dealII91}.

\begin{figure}[ht]
\centering
~\hfill{
	\includegraphics[width=4cm]{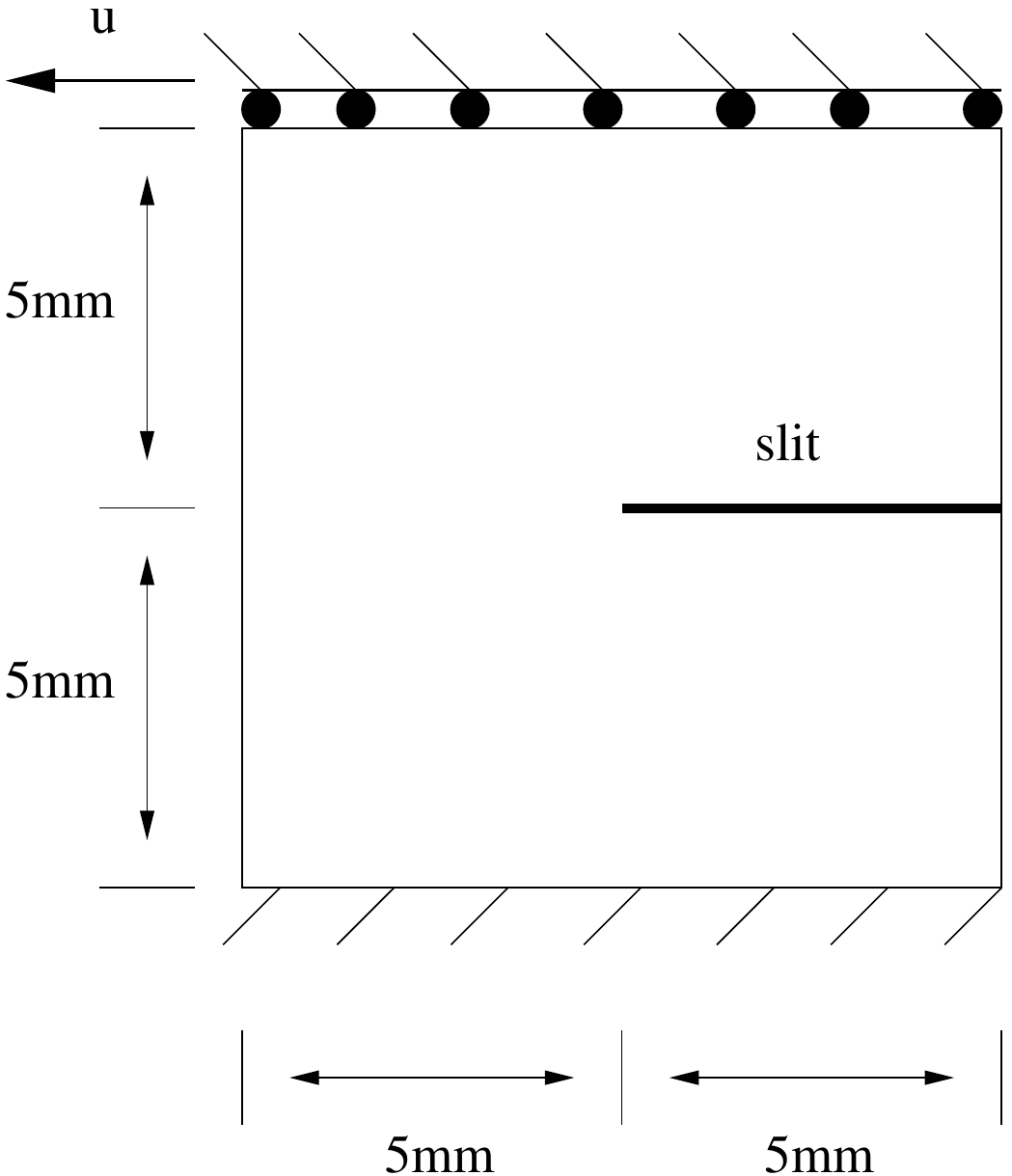}}\hfill\hfill
{\includegraphics[width=5cm]{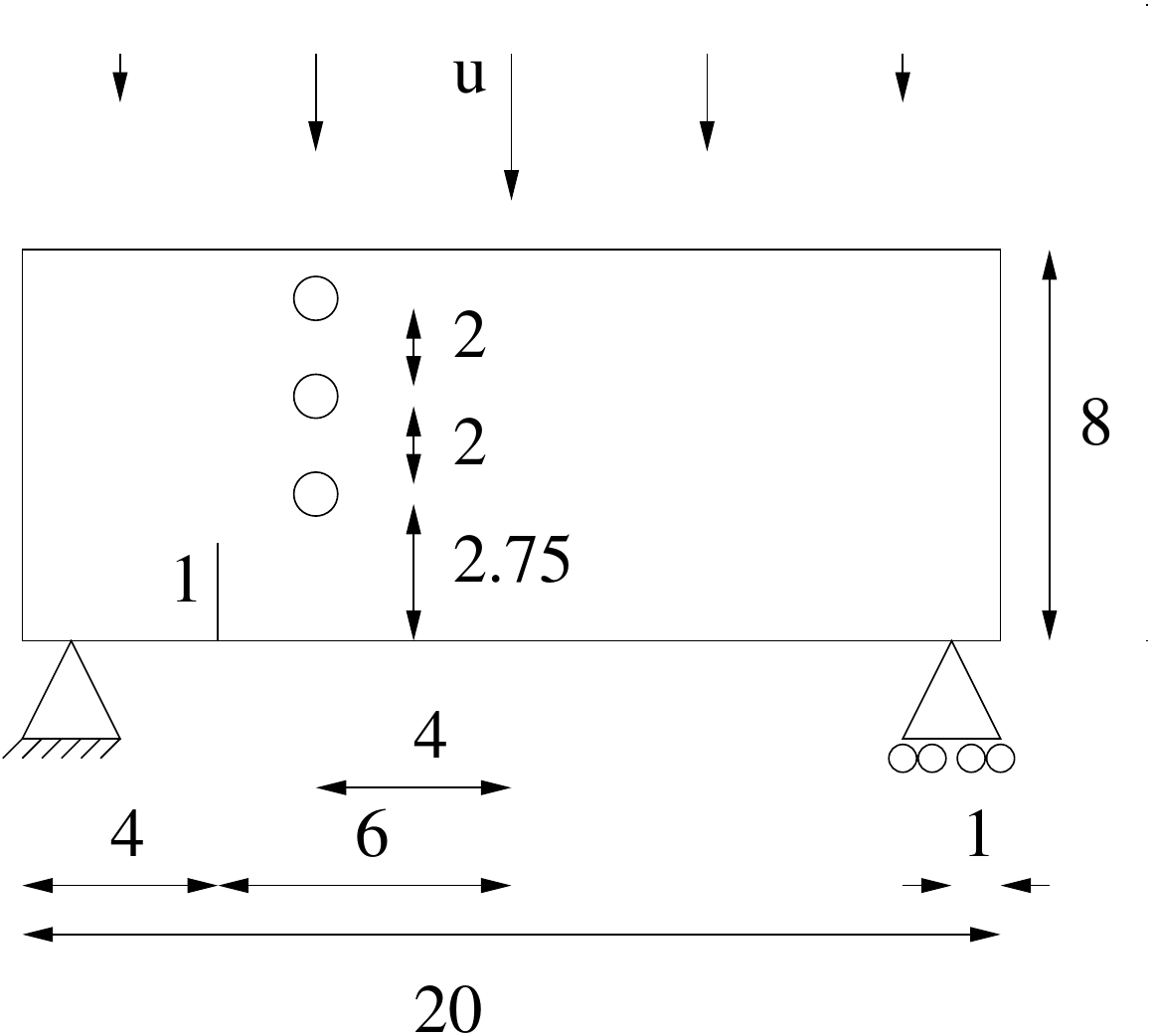}}
\hfill~
\caption{Examples 1 and 2. %Configurations.
The following conditions are prescribed: on the left and right
boundaries, $u_y$ = \SI{0}{mm} and traction-free in $x$-direction. 
On the bottom 
part, $u_x = u_y$ = \SI{0}{mm}. On $\Gamma_\textnormal{top}$, $u_y$ = \SI{0}{mm} 
and $u_x$ is as stated in \eqref{diri_ex_2}. Finally, the lower part of the 
slit is fixed in $y$-direction, i.e., $u_y$ = \SI{0}{mm}. 
Right: Asymmetric notched three point bending test. 
The three holes have 
each a diameter of $0.5$. All units are in ${mm}$. 
}
\label{ex_1_config}
\end{figure}

\begin{figure}[h!]
\centering
~\hfill{\includegraphics[width=4cm]{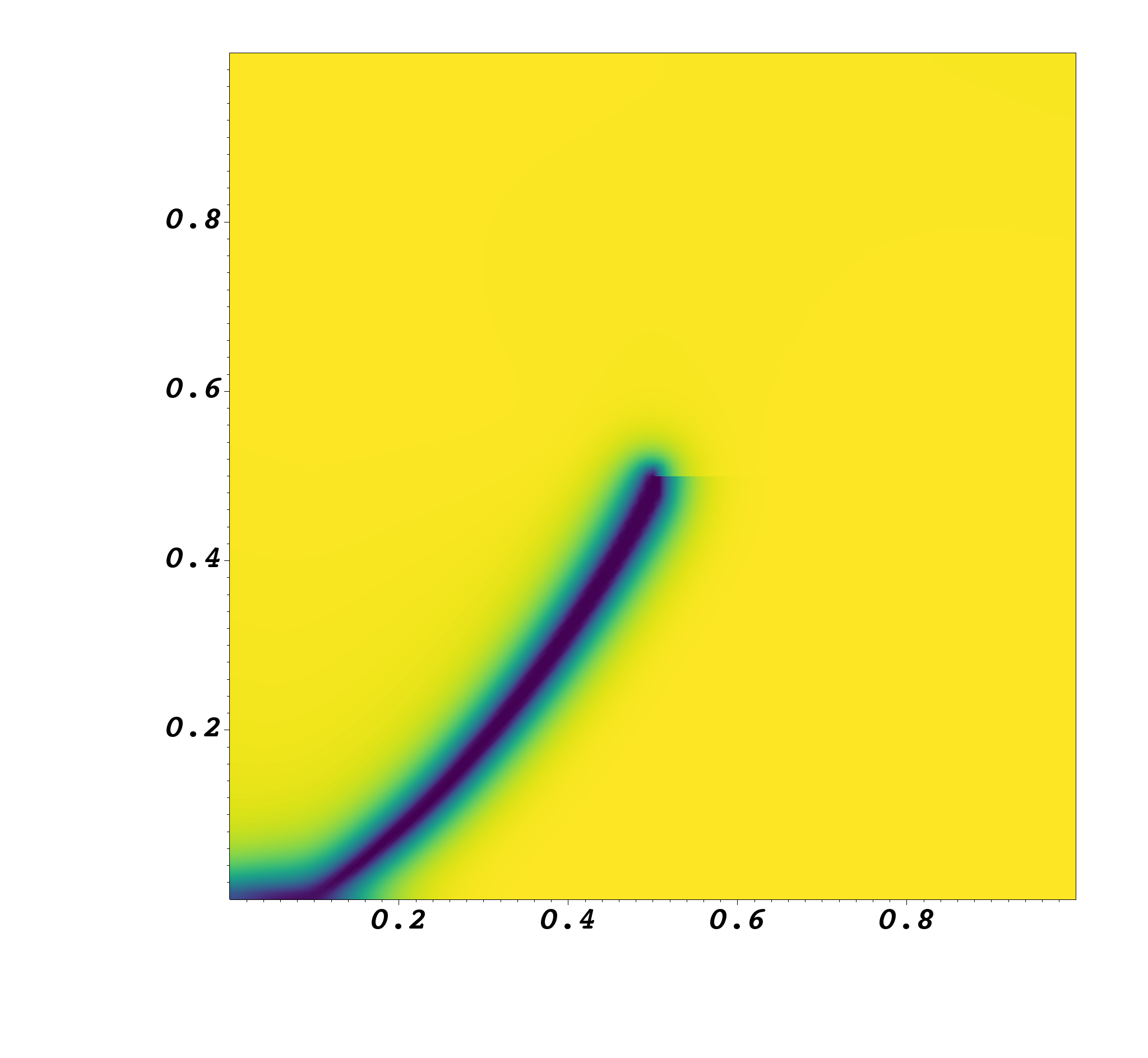}}\hfill\hfill
{\includegraphics[width=6cm]{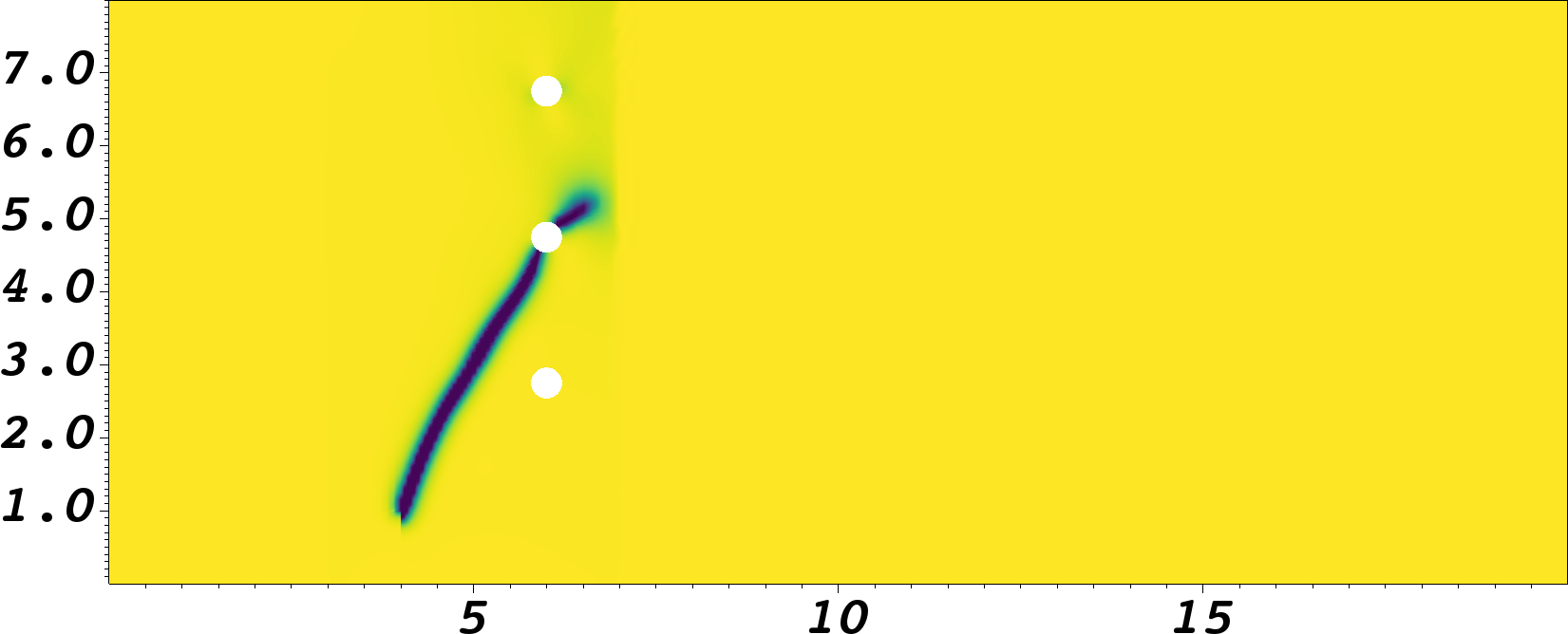}}\hfill~
\caption{Examples 1 and 2. Numerical solutions on the finest meshes
and at the end time. The cracks are displayed in dark blue color.
}
\label{ex_1_num_sol}
\end{figure}

\medskip\noindent\textbf{Single edge notched shear test:\quad}
The configuration is shown in Figure \ref{ex_1_config}.
Specifically, we use
$\mu_s$ = \SI{80.77}{kN/mm^2}, $\lambda_s$ = \SI{121.15}{kN/mm^2}, and
$G_c$ = \SI{2.7}{N/mm}.
The crack growth is driven  by a non-homogeneous Dirichlet
condition for the displacement field on $\Gamma_\textnormal{top}$, the top boundary of $B$ .
We increase the displacement on $\Gamma_\textnormal{top}$ over time, namely we apply
non-homogeneous Dirichlet conditions:
\begin{align}
u_x &= t \bar{u}, \quad \bar{u} = \SI{1}{mm/s}, \label{diri_ex_2}
\end{align}
where $t$ denotes the current loading time. Furthermore, we set $\kappa = 10^{-10}$ [mm] and $\eps = 2h$ [mm]. 
We evaluate the surface load vector on the $\Gamma_\textnormal{top}$ as
\begin{equation}
\label{eq_Fx_Fy}
\tau = (F_x,F_y) := \int_{\Gamma_\textnormal{top}} \bsig(u)\nu\, \textnormal{d}s,
\end{equation}
with normal vector $\nu$,
and we are particularly interested in the shear force $F_x$. 
Three different meshes with $1024$ (Ref. 4), $4096$ (Ref. 5) and $16384$ (Ref. 6) 
elements are observed 
in order to show the robustness of the proposed schemes. The results 
are shown in Figure \ref{ex_1_and_2_mesh_levels}.

Our findings are summarized in Figure \ref{ex_1_b}.
The numerical solutions for all four different strategies for choosing $L$
are practically identical, only the number of iterations being different. Here,  
$L=0$ and $L=1e-2$ denote tests in which $L=L_u=L_{\phi}$ are taken constant
throughout the entire computation. The newly proposed dynamic
versions are denoted by \texttt{L dynamic} and \texttt{L dyn. weighted}.
We observe a significant reduction in the computational 
cost when using the dynamic $L$-schemes. The maximum number of 
iterations is $21$ for both the weighted version 
and the spatially-constant $L$-scheme.
This number is reduced to $12$ iterations using $a=20$ while 
the accuracy only slightly changes.

\begin{figure}[H]
\centering
{\includegraphics[width=5cm]{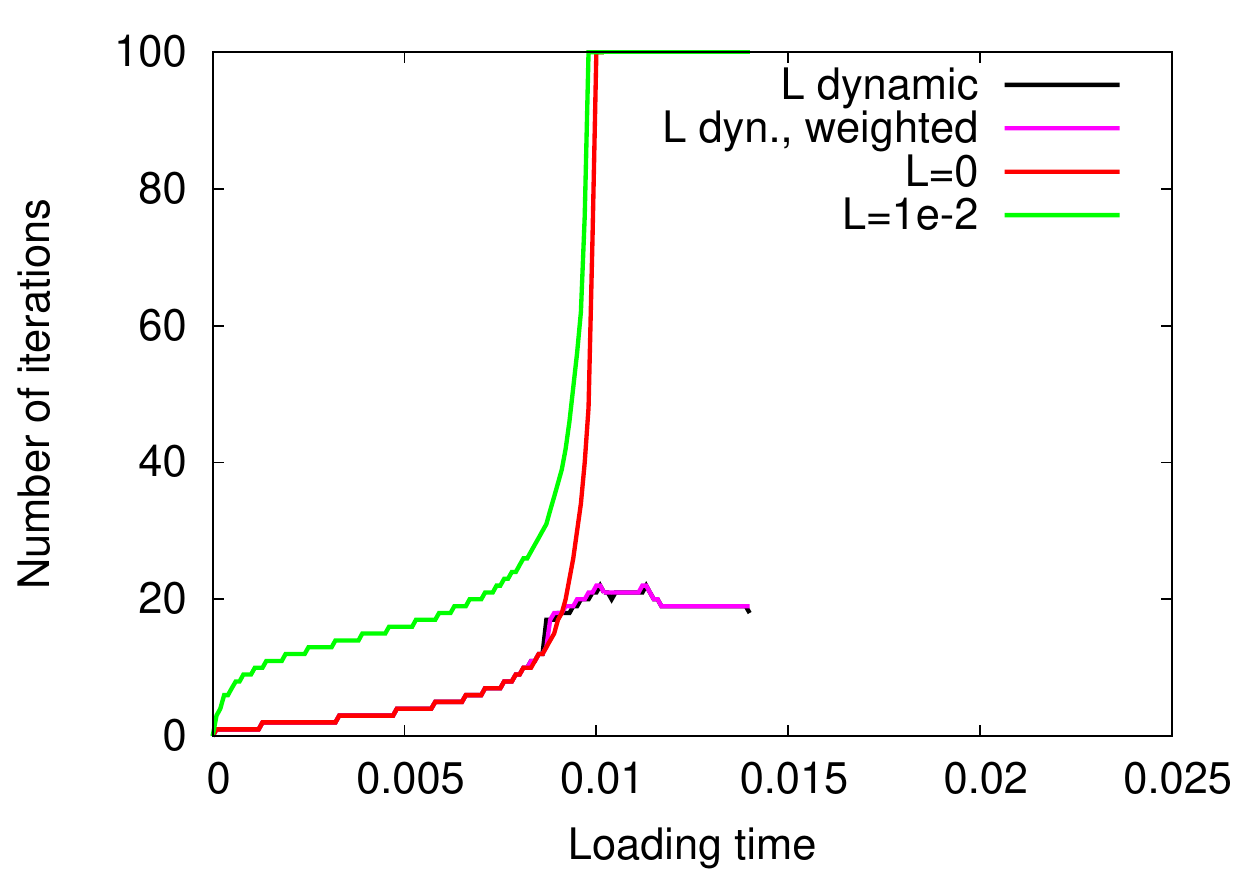}} 
{\includegraphics[width=5cm]{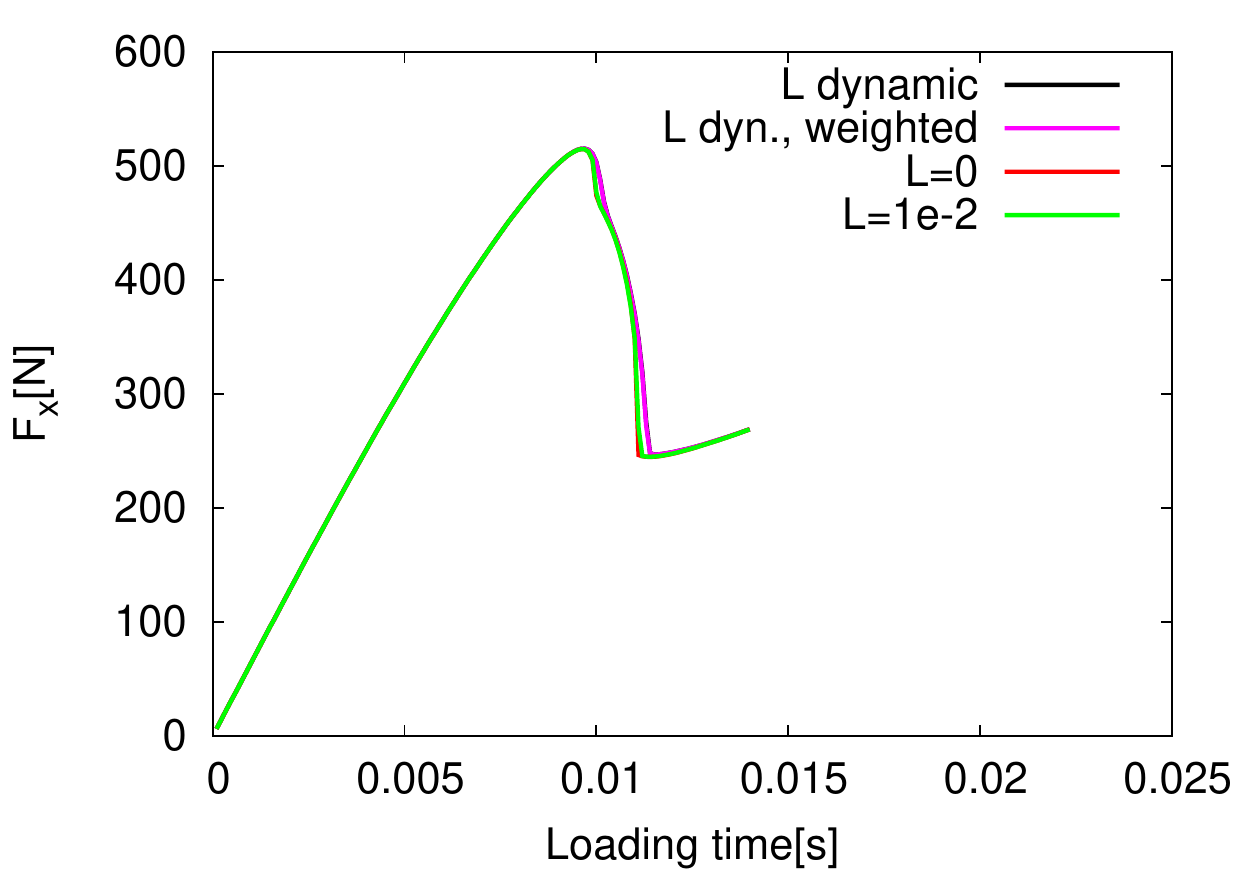}} 
%\added[id=sp]{Sorry, ich habe die Bilder nicht}
\caption{Example 1. Comparison of dynamic $L$ updates, the weighted version,
  and constant $L$.
  Left:~number of iterations.
  Right:~load-displacement curves.}
\label{ex_1_b}
\end{figure}

\begin{figure}[H]
\centering
{\includegraphics[width=5cm]{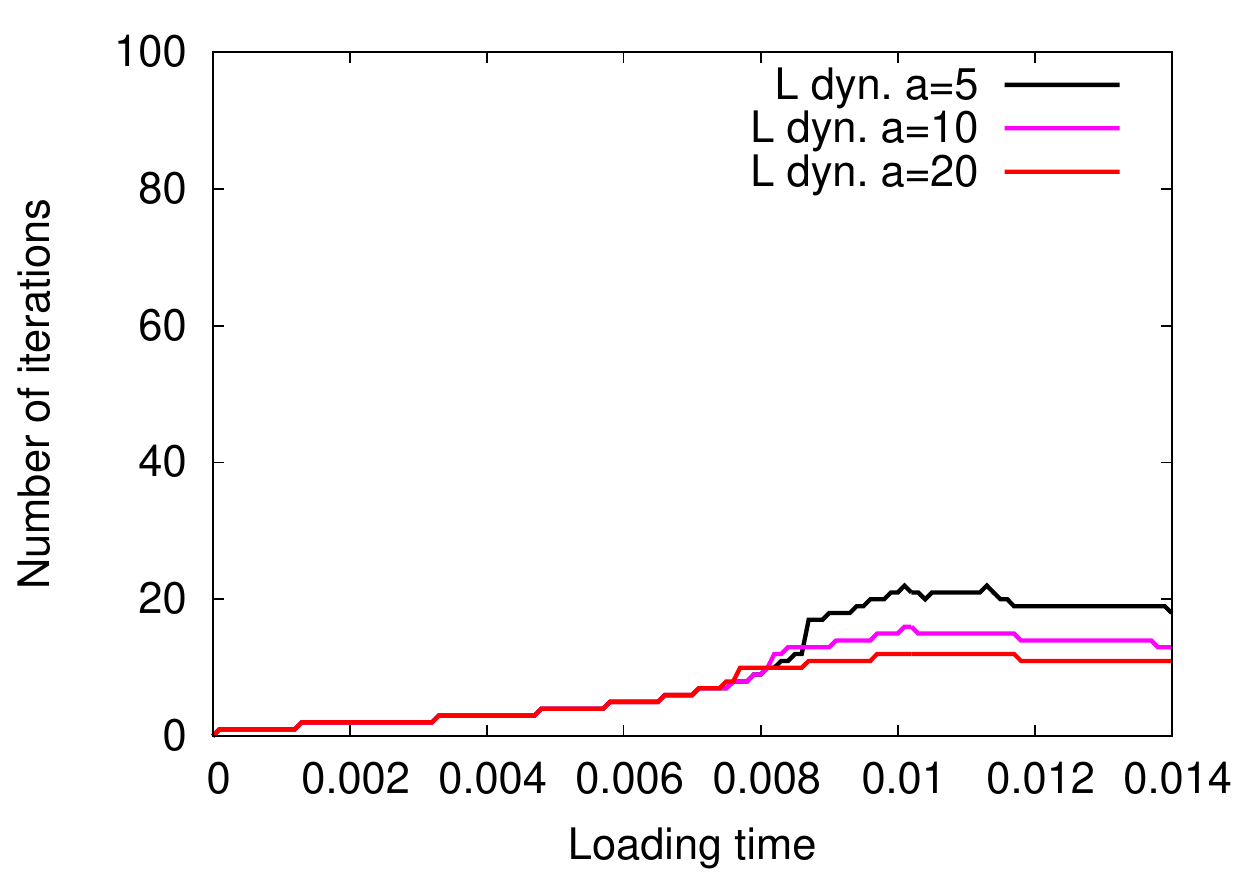}} 
{\includegraphics[width=5cm]{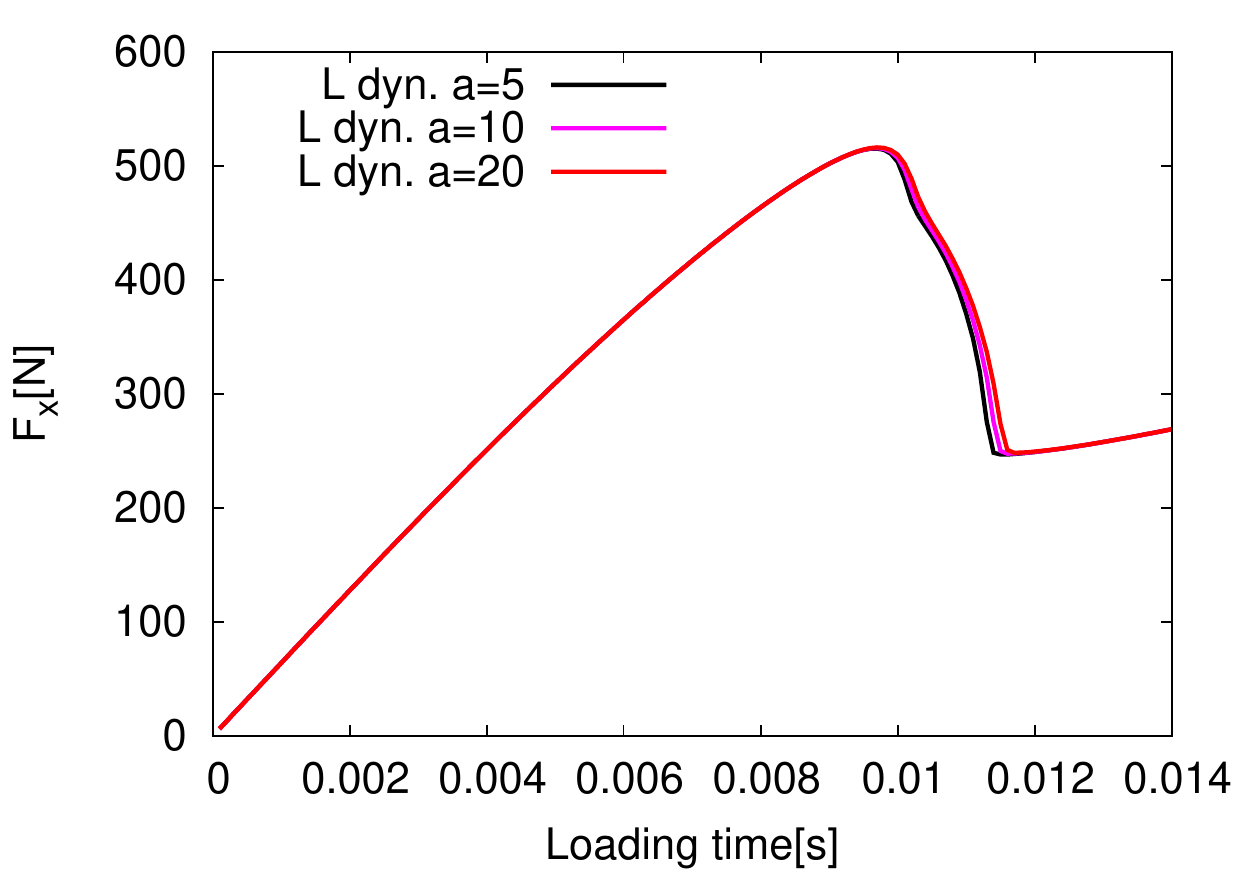}} 
\caption{Example 1. Comparison of different $a$ for the dynamic $L$ scheme.}
\label{ex_1_c}
\end{figure}

\newpage
\medskip\noindent\textbf{Asymmetrically notched three point bending test:\quad}
The configuration is shown in Figure \ref{ex_1_config} (right).
The initial mesh is $3,4$ and $5$ times uniformly refined, yielding 
$3\, 904, 15\, 616$ and $62\, 464$ mesh elements 
with the minimal mesh size parameter 
$h_3=0.135, h_4 = 0.066$ and $h_5 = 0.033$.
As material parameters, we use
$\mu_s$ = \SI{8}{kN/mm^2}, $\lambda_s$ = \SI{12}{kN/mm^2}, and
$G_c$ = \SI{1e-3}{kN/mm}. 
Furthermore, we set $k = 10^{-10}h$[mm] and 
$\eps=2h$. 

%\begin{figure}[H]
%\centering
%{\includegraphics[width=6cm]{pic_asymmetric.pdf}} 
%%\hspace*{0.5cm}
%%{\includegraphics[width=9cm]{visit_Sep_18_2019_0000.png}} 
%\caption{Example 2: Asymmetric notched three point bending test. 
%The three holes have 
%each a diameter of $0.5$. All units are in ${mm}$. At right, the three times 
%uniformly refined mesh is displayed.}
%\label{asymmetric_a}
%\end{figure}

Figure \ref{asymmetric_b} presents the number 
of iterations and the load-displacement curves. The number of iterations is decreasing from $500$ (in the figures cut to $100$) 
for the classical L-scheme, to a maximum of $25$ when using the dynamic updates. 
The choice of weighting does not seem to have 
a {significant} influence on the number of iterations though.
The crack starts growing a bit later when using 
the dynamic updates, which can be inferred from 
the {right plot} in Figure \ref{asymmetric_b}. Thus, the 
stabilization parameters have a slight influence on the 
physical solution. This can be explained in the following way. In regions where $\phi=0$ the
solution component $u$ is not uniquely defined. This leads to a sub-optimal convergence behaviour of the L-scheme. 
With the dynamic L-scheme we regain uniqueness, but at the cost of a slightly modified physical problem.

\begin{figure}[H]
\centering
{\includegraphics[width=5cm]{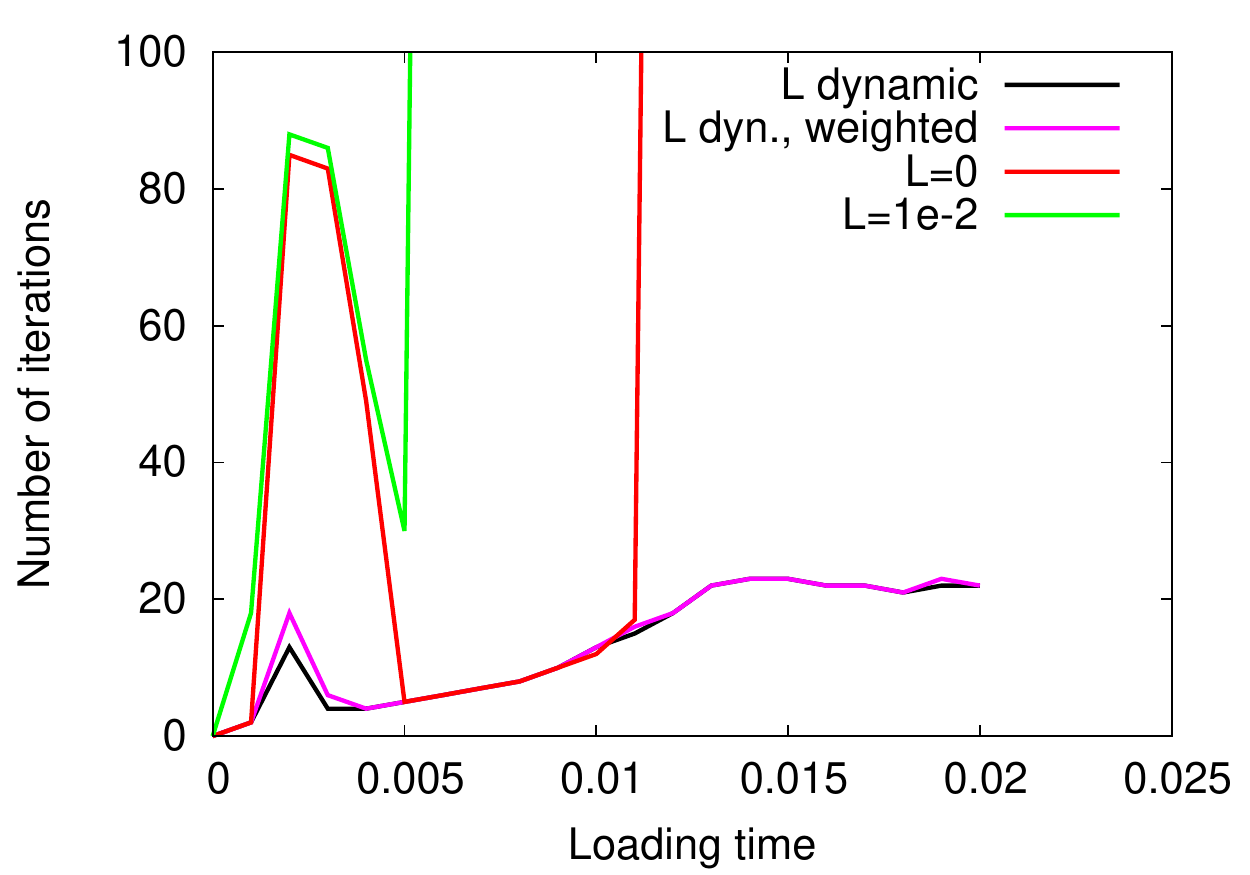}} 
{\includegraphics[width=5cm]{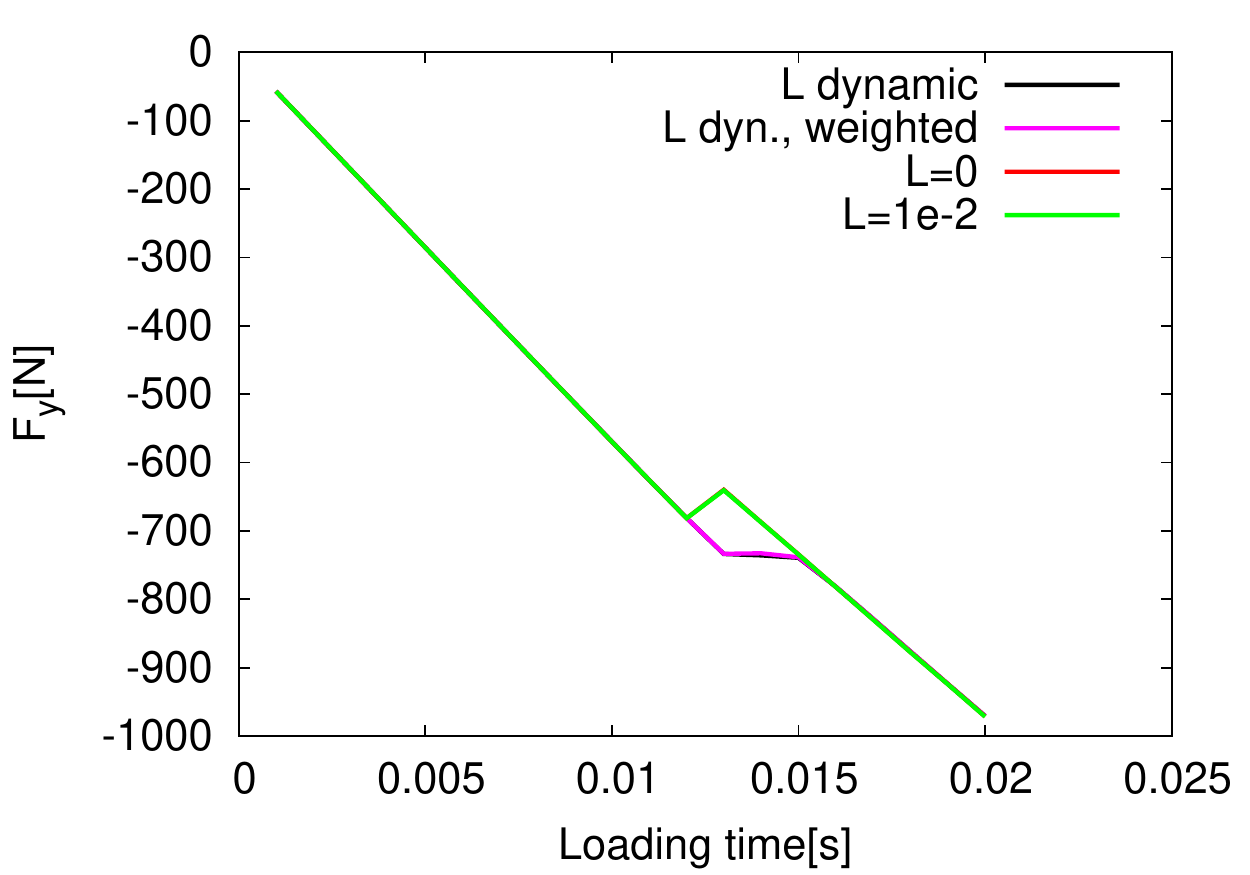}} 
%\added[id=sp]{Sorry, ich habe die Bilder nicht}
\caption{Example 2: Left: The number of iterations for the different schemes; 
	the results for $L=0$ and $L=1e-2$ are taken from \cite{BrWiBeNoRa19}.
	Right: The load-displacement curves; a slight difference can be observed in the results, 
indicating that the dynamic updates lead to a slight delay in the prediction 
of the starting time for the fracture growth.
}
\label{asymmetric_b}
\end{figure}

\begin{figure}[H]
\centering
{\includegraphics[width=5cm]{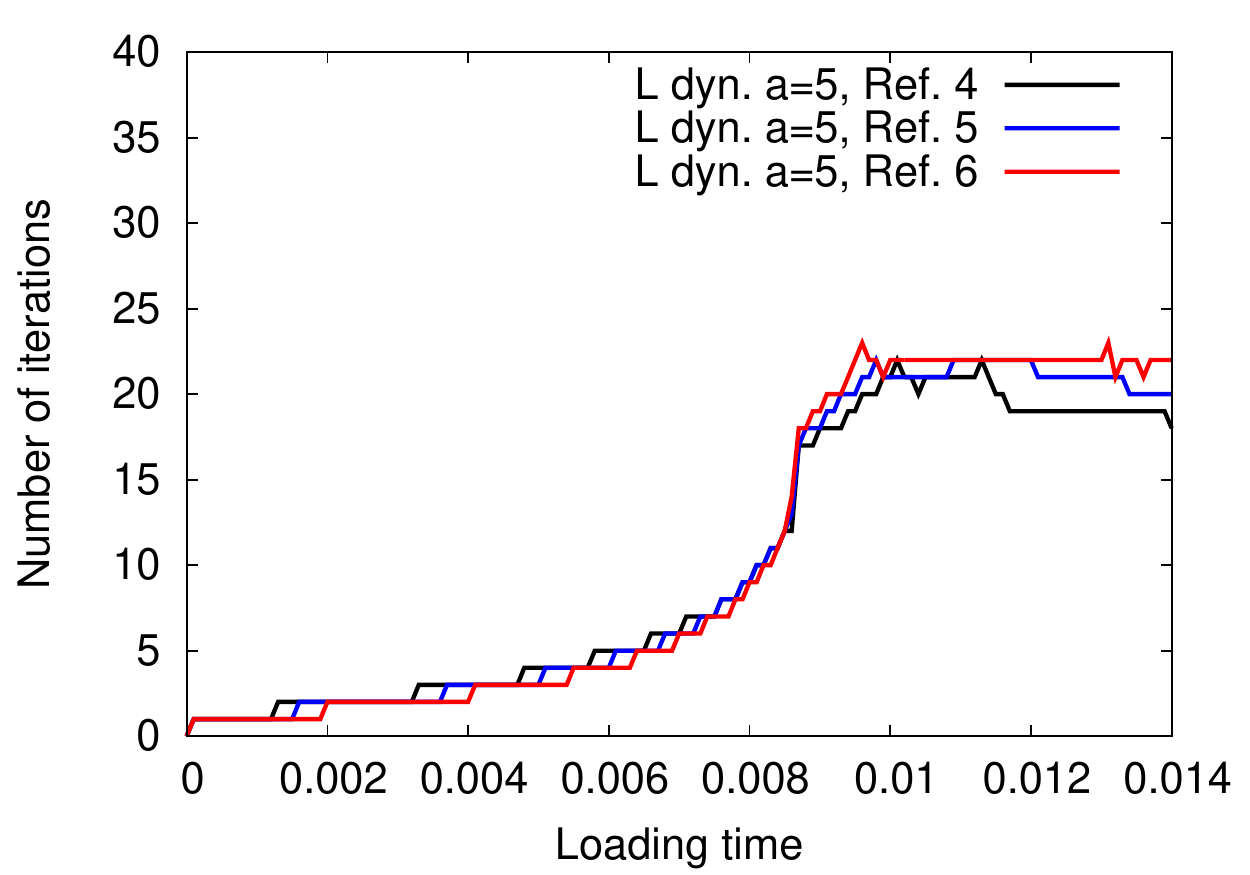}} 
{\includegraphics[width=5cm]{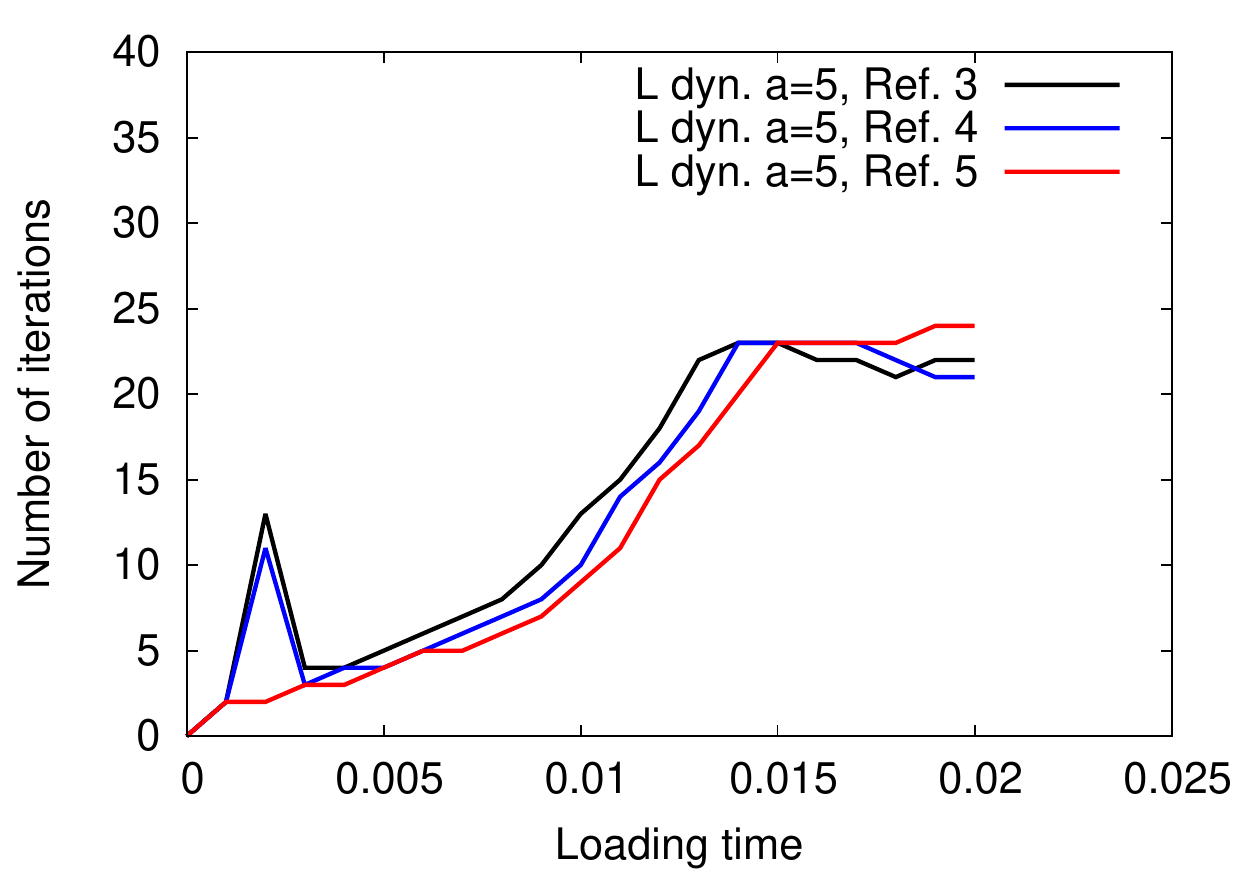}} 
%\added[id=sp]{Sorry, ich habe die Bilder nicht}
\caption{Examples 1 and 2 for the dynamic $L$ scheme using $a=5$; three 
different mesh levels are used in order to verify the 
robustness of the proposed scheme. The results indicate 
that the mesh size does not influence the number of the iterations.}
\label{ex_1_and_2_mesh_levels}
\end{figure}

%\begin{remark}
%\added[id=ce]{
%This can be interpreted in the following way. Where $\phi=0$ the
%solution of $u$ is not uniquely defined. This leads to the observed
%issues regarding the convergence. With our dynamic L-scheme we
%regain uniqueness, but at the cost of a slightly modified physical problem.
%}
%\end{remark}

\begin{remark}
Noteworthy, the number of iterations for the
  dynamic L-scheme is robust
  with respect to the mesh refinement, as
shown in Figure 
\ref{ex_1_and_2_mesh_levels}. This is in line with the analysis in \cite{BrWiBeNoRa19, list2016study, MR2079503}, where it is proved that the convergence rate does not depend on the spatial discretization. 
\end{remark}

%%%%%%%%%%%%%%%%%%%%%%%%%%%%%%%%%%%%%%%%%%%%%%%%%%%%%%%%%%
\section*{Acknowledgements}\small
TW is supported by the German Research Foundation, Priority Program 1748 (DFG SPP 1748) 
under the project No.~392587580.
CE is supported by the German Research Foundation, via
  Priority Program 1648  (DFG SPP 1648) under the grant
  No.~EN-1042/2-2 and via EXC 2044-390685587, Mathematics
  M{\"u}nster: Dynamics-Geometry-Structure.
ISP is supported by the Research Foundation-Flanders (FWO), Belgium through the Odysseus programme (project G0G1316N).

%%%%%%%%%%%%%%%%%%%%%%%%%%%%%%%%%%%%%%%%%%%%%%%%%%%%%%%%%%
%\input{references}
\bibliographystyle{abbrv}
\bibliography{lit}

\end{document}